# A Two-Point Newton Method suitable for non-convergent Cases and with Super-Quadratic Convergence

Ababu Teklemariam Tiruneh[1]

[1]Lecturer, Department of Environmental Health Science. University of Swaziland.
P.O.Box 369, Mbabane H100, Swaziland.  Email: ababute@yahoo.com

**Abstract**

An iterative formula based on Newton's Method alone is presented for the iterative solutions of equations that ensures convergence in cases where the traditional Newton Method may fail to converge to the desired root. In addition, the method has super quadratic convergence of order 2.414 (i.e., 1+ $\sqrt{2}$ ). Newton method is said to fail in certain cases leading to oscillation, divergence to increasingly large number or off-shooting away to another root further from the desired domain or off shooting to an invalid domain where the function may not be defined. In addition when the derivative at the iteration point is zero, Newton method stalls. In most of these cases, hybrids of several methods such as Newton, bisection and secant methods are suggested as substitute methods and Newton method is essentially blended with other methods or altogether abandoned. This paper argues that a solution is still possible in most of these cases by the application of Newton Method alone without resorting to other methods and with the same computational effort (two functional evaluations per iteration) like the traditional Newton method.  In addition, the proposed modified formula based on Newton method has better convergence characteristics than the traditional Newton method.

## 1. Introduction

Iterative procedures for solutions of equations are routinely employed in many science and engineering problems.  Starting with the classical Newton methods, a number of methods for finding roots of equations have come to exist each of which has its own advantages and limitations. The Newton method of root finding is based on the iterative formula:

$$x_{k+1} = x_k - \frac{y(x_k)}{y'(x_k)}$$

Newton's method displays a faster quadratic convergence near the root while it requires evaluation of the function and its derivative at each step of the iteration. However, when the derivative evaluated is zero, Newton method stalls. For low values of the derivative, the Newton iteration offshoots away from the current point of iteration and may possible converge to a root far away from the intended domain. For certain forms of equations, Newton method diverges or oscillates and fails to converge to the desired root. In addition, the convergence of Newton method can be slow near roots of multiplicity although modifications can be made to increase the rate of convergence [1].

Modification of the Newton method with higher order convergence have been proposed that require also evaluation of a function and its derivatives. An example of such methods is a third order convergence method by S. Weeraksoon and T.G. Fernando [2] that requires evaluation of one function and two first derivatives. A fourth order iterative method, accordintg to J.F. Traub [3] also requires evaluation of one function and two derivative evaluations. Sanchez and Barrero [4] gave a compositing of function evaluation at a point and its derivative to improve the convergence of Newton's method from 2 to 4. Recently other methods of fifth, sixth, seventh and higher order convergence have been proposed [5-11]. While higher order Newton methods ensure faster convergence, their stability for certain equation forms may have the same problem similar to that of the taditional Newton method.

The secant method does not require evaluation of derivatives. However, the rate of convergence is about 1.618 and the convergence may be a problem for some forms of equations and the secant method may fail to converge in those cases. Muller's method is an extension of the secant method to a quadratic polynomial [12]. It requires three functional evaluations to start with but continues with one function evaluation afterwards. The method does not require derivatives and the rate of convergence is about 1.84. However, Muller's method can converge to a complex root from an initial real number [13]. Muller's method also requires the three points to be distinct. If two of the points coincide, the method degenerates to secant method.

Hybrids of methods are also used to provide stability and ensure convergence to a desired root. For example Newton or Secant methods can be combined with bisection to bracket the root by a small interval so that a good initial guess is available for applying Newton method with quadratic convergence. Dekker's method [14] combines bisection method with that of secant method. The method starts by bracketing the root between two initial points that have functional values opposite in sign. The secant estimate of x is compared with the bisection of the interval and the one estimate resulting in small interval with the point of smaller magnitude functional value is chosen and the iteration continues. Brent's method [15] is a root finding algorithm that combines root bracketing, bisection and inverse quadratic interpolation. It is a modification of Dekker's method to avoid slow convergence when the difference between consecutive estimates of x is arbitrarily small. In such cases bisection is used for the next root estimate.

The Leap-frogging Newton method [16] uses the Newton method as an intermediate step for the estimation of the root followed by the application of the secant method using the initial point and the intermediate point found by Newton method. The method has cubic convergence and works in certain pathological cases where Newton method may fail. However, being dependent on Newton method for the intermediate step, the method may suffer from the same drawbacks of using the traditional Newton method mentioned above.

## 2. Method development

It will be shown that the iterative formula for the two point method will take the form;

$$x_{k+1} = x_{k-1} - \frac{(x_{k-1} - x_k)}{1 - \left(\frac{y_k}{y_{k-1}}\right)\left(\frac{\frac{y_k - y_{k-1}}{x_k - x_{k-1}}}{y'_k}\right)} \qquad (1)$$

The method starts by selecting two points lying on a curve y= f(x), namely ($x_0$, $y_0$) and ($x_1$, $y_1$). A line connecting the two points is drawn from the point ($x_0$, $y_0$) to the point ($x_1$, $y_1$) as shown in Figure 1 below. A new variable m (the cotangent of the angle $\alpha_1$ between this line and the vertical) is defined so that:

$$m_1 = \cot(\alpha_1) = \frac{y_0 - y_1}{x_0 - x_1}$$

In general, for any other point (x, y) lying on the curve, the variable m is defined so that:

$$m = \cot(\alpha) = \frac{y_0 - y}{x_0 - x}$$

Next, Newton Method of finding roots will be applied taking m as the independent variable and y as the dependent variable. Denoting that $m_r$ is the estimate of the root for which y (m) = 0 from the Newton method with the corresponding value of $x_r$ and applying Newton method gives;

$$m_r = m_1 - \frac{y_1}{\left(\frac{dy}{dm_1}\right)} \qquad (2)$$

Substituting for $m_r$ and $m_1$ the equivalent expressions in terms of x and y values, namely,

$$m_r = \frac{y_0 - 0}{x_0 - x_r} \qquad \text{and} \qquad m_1 = \frac{y_0 - y_1}{x_0 - x_1}$$

into Equation 2 above gives the following expression;

$$\frac{y_0}{x_0 - x_r} = \left(\frac{y_0 - y_1}{x_0 - x_1}\right) - \frac{y_1}{\left(\frac{dy}{dm_1}\right)} \qquad (3)$$

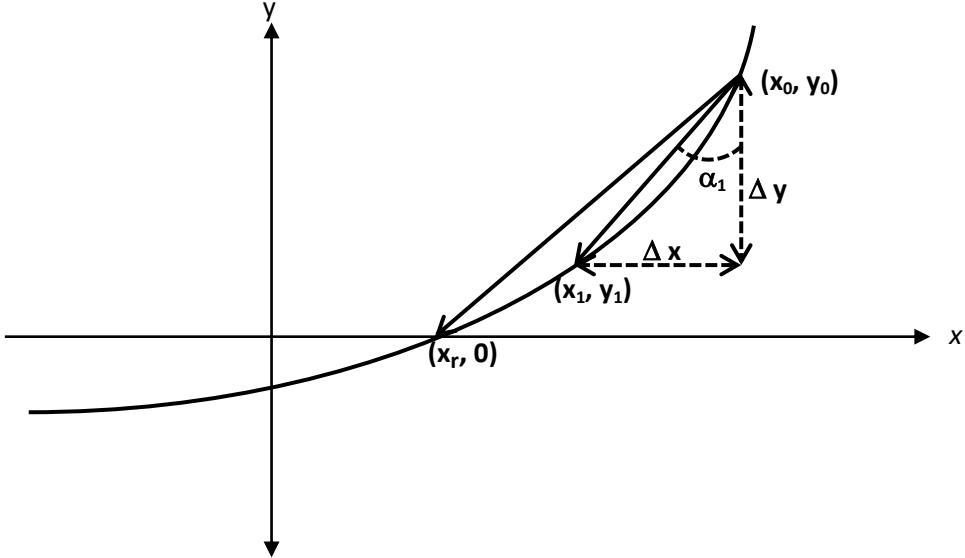

Figure 1  The graph of y = f(x) on which the angle α is defined between the two points.

The derivative dy/dm is evaluated from dy/dx = y' and dm/dx using the formula;

$$\frac{dy}{dm} = \frac{dy/dx}{dm/dx} = \frac{y'}{dm/dx}$$

The evaluation of the derivative dm/dx in turn gives;

$$\frac{dm}{dx} = \frac{d}{dx}\left(\frac{y_0 - y}{x_0 - x}\right) = \frac{(y_0 - y) - (x_0 - x)y'}{(x_0 - x)^2}$$

Therefore, dy/dm can now be written in terms of x, y and y' as follows;

$$\frac{dy}{dm} = \frac{y'(x_0 - x)^2}{(y_0 - y) - (x_0 - x)y'}$$

Substituting the above expression for dy/dm in the equation,

$$\frac{y_0}{x_0 - x_r} = \left(\frac{y_0 - y_1}{x_0 - x_1}\right) - \frac{y_1}{\left(\frac{dy}{dm_1}\right)}$$

gives;

$$\frac{y_0}{x_0 - x_r} = \left(\frac{y_0 - y_1}{x_0 - x_1}\right) - \left[\frac{y_1}{\left(\frac{y'(x_0 - x)^2}{(y_0 - y) - (x_0 - x)y'}\right)}\right]$$

Solving for the root estimate $x_r$ and further rearranging results in;

$$x_r = x_0 - \frac{(x_0 - x_1)}{1 - \left(\frac{y_0}{y_1}\right)\left(\frac{\frac{y_1 - y_0}{x_1 - x_0}}{y'_1}\right)}$$

Continuing the iteration using the above formula, for the $k^{th}$ step of the iteration, the $k+1^{th}$ estimate of the root will take the form;

$$x_{k+1} = x_{k-1} - \frac{(x_{k-1} - x_k)}{1 - \left(\frac{y_k}{y_{k-1}}\right)\left(\frac{\frac{y_k - y_{k-1}}{x_k - x_{k-1}}}{y'_k}\right)} \qquad (4)$$

Denoting $(y_k - y_{k-1})/(x_k - x_{k-1})$ as $(\Delta y/\Delta x)_k$ will result in the expression;

$$x_{k+1} = x_{k-1} - \frac{(x_{k-1} - x_k)}{1 - \left(\frac{y_k}{y_{k-1}}\right)\left(\frac{(\frac{\Delta y}{\Delta x})_k}{y'_k}\right)}$$

## 3. Proof of super-quadratic convergence

Recalling the iteration formula of Equation 4 above:

$$x_{k+1} = x_{k-1} - \frac{(x_{k-1} - x_k)}{1 - \left(\frac{y_k}{y_{k-1}}\right)\left(\frac{\frac{y_k - y_{k-1}}{x_k - x_{k-1}}}{y'_k}\right)}$$

Defining the error at the $k^{th}$ iteration to be $E_k = X_k - r$ where r is the root of the equation desired, the errors at the $k-1^{th}$ and $k^{th}$ iteration are also defined similarly:

$$E_{k-1} = x_{k-1} - r \text{ and } E_{k+1} = x_{k+1} - r$$

The iteration formula in equation 4 above can now be rewritten in terms of the error terms as follows:

$$E_{k+1} = E_{k-1} - \frac{(E_{k-1} - E_k)}{1 - \left(\frac{y_k}{y_{k-1}}\right)\left(\frac{\frac{y_k - y_{k-1}}{E_k - E_{k-1}}}{y'_k}\right)} \qquad (5)$$

Expanding the $y_{k-1}$, $y_k$ and $y'_k$ terms about the root r using Taylor Series expansion;

$$y_k = (y_r = 0) + E_k y'_r + \frac{E_k{}^2 y''_r}{2} + O(E_k)^3$$

$$y_{k-1} = (y_r = 0) + E_{k-1} y'_r + \frac{E_{k-1}{}^2 y''_r}{2} + O(E_k)^3$$

$$y'_k = y'_r + E_k y''_r + \frac{E_k{}^2 y'''_r}{2} + O(E_k)^3$$

$$\frac{(y_k - y_{k-1})}{E_k - E_{k-1}} = y'_r + \left(\frac{E_k + E_{k-1}}{2}\right) y''_r + O(E_k)^3$$

$$\frac{y_k}{y_{k-1}} = \frac{(y_r = 0) + E_k y'_r + \frac{E_k{}^2 y''_r}{2} + O(E_k)^3}{(y_r = 0) + E_{k-1} y'_r + \frac{E_{k-1}{}^2 y''_r}{2} + O(E_k)^3}$$

and

$$1 - \left(\frac{y_k}{y_{k-1}}\right)\left(\frac{\frac{y_k - y_{k-1}}{E_k - E_{k-1}}}{y'_k}\right) = 1 - \left(\frac{(y_r = 0) + E_k y'_r + \frac{E_k^2 y''_r}{2} + O(E_k)^3}{(y_r = 0) + E_{k-1} y'_r + \frac{E_{k-1}^2 y''_r}{2} + O(E_k)^3}\right)\left(\frac{y'_r + \left(\frac{E_k + E_{k-1}}{2}\right) y''_r + O(E_k)^3}{y'_r + E_k y''_r + \frac{E_k^2 y'''_r}{2} + O(E_k)^3}\right)$$

................................ (6)

Defining $c = E_{k-1}$ and assuming that the convergence is of order n so that $E_k = E_{k-1}^n$ and $E_{k+1} = E_k^n$, etc., both $E_k$ and $E_{k+1}$ can now be written in terms of $E_{k-1} = c$;

$$E_{k-1} = c \; ; \; E_k = c^n \text{ and } E_{k+1} = c^{n^2}$$

$$1 - \left(\frac{y_k}{y_{k-1}}\right)\left(\frac{\frac{y_k - y_{k-1}}{E_k - E_{k-1}}}{y'_k}\right)$$

$$= \frac{cy'_r{}^2 + \frac{c^2 y''_r y'_r}{2} + c^{n+1} y'_r y''_r + \frac{c^{2+n} y''_r{}^2}{2} + \frac{c^{2n+1} y'_r y''_r}{2} + \frac{c^{2n+2} y''_r y'''_r}{4} - c^n y'_r{}^2 - c^{2n} y'_r y''_r - \frac{c^{n+1} y'_r y''_r}{2} - \frac{c^{2n+1} y''_r{}^2}{4}}{cy'_r{}^2 + \frac{c^2 y''_r y'_r}{2} + c^{n+1} y''_r y'_r + \frac{c^{n+2} y''_r{}^2}{2} + \frac{c^{2n+1} y'_r y''_r}{2} + \frac{c^{2n+2} y''_r y'''_r}{4}}$$

................................ (7)

Substituting the right hand side of the above equation (Equation 7) in the iteration formula of Equation 5 above and simplifying one gets the following expression:

$$E_{k+1} = c^{n^2}$$

$$= \frac{\frac{c^{2n+2} y''_r{}^2}{4} + \frac{c^{3n+1} y'_r y'''_r}{2} + \frac{c^{3n+2} y''_r y'''_r}{4}}{cy'_r{}^2 + \frac{c^2 y''_r y'_r}{2} + c^{n+1} y'_r y''_r + \frac{c^{2+n} y''_r{}^2}{2} + \frac{c^{2n+1} y'_r y''_r}{2} + \frac{c^{2n+2} y''_r y'''_r}{4} - c^n y'_r{}^2 - c^{2n} y'_r y''_r - \frac{c^{n+1} y'_r y''_r}{2} - \frac{c^{2n+1} y''_r{}^2}{4}}$$

Dividing all the expressions by c will leave the denominator with the dominant expression $(y'_r)^2$ since all the other expression in the denominator contain the error term $c = E_{k-1}$. Discarding the terms containing c in the denominator after division by c and further simplification will give:

$$c^{n^2} = \frac{c^{2n+1} \left(\frac{y''_r}{y'_r}\right)^2}{4} + \frac{c^{3n} y''_r}{2 y'_r} + \frac{c^{3n+1} y''_r y'''_r}{4 y'_r{}^2}$$

Equating the greatest error term in the right side, $c^{2n+1}$, with that of $c^{n^2}$ will give the equation:

$$n^2 - 2n - 1 = 0$$

$$n = \frac{2 + \sqrt{4 + 4}}{2} = 1 + \sqrt{2} \approx 2.414$$

Therefore, the iteration formula of Equation 1 near the root has a convergence of order 2.414 which is super quadratic.

## 4. Stability of the proposed method

The proposed method shows stability near points where application of the traditional Newton Method may result in oscillation, divergence or off-shooting away from the desired root or off-shooting to a possibly invalid domain where the function y=f(x) may not be defined. This stability characteristics may be studied by examining the iteration formula again which takes the following form;

$$x_{k+1} = \left(1 - \frac{1}{r}\right) x_{k-1} + \left(\frac{1}{r}\right) x_k$$

Where the variable r is defined as;

$$r = 1 - \left(\frac{y_k}{y_{k-1}}\right)\left(\frac{\frac{y_k - y_{k-1}}{x_k - x_{k-1}}}{y'_k}\right)$$

The estimate for the root at the k+1 $^{th}$ iteration, $x_{k+1}$, therefore, can be taken as the weighted sum of the $x_{k-1}$ and $x_k$ values at the k-1$^{th}$ and k$^{th}$ iteration respectively. The weighing factors are (1-1/r) for $x_{k-1}$ and (1/r) for $x_k$. When the method converges to the root the value of r approaches unity. This is found by examining the behavior of r near the root and noting that $y_k$ approaches 0 near the root, i.e.,

$$\lim_{k \to \infty}(r) = 1 - \left(\frac{0}{y_{k-1}}\right)\left(\frac{\frac{y_k - y_{k-1}}{x_k - x_{k-1}}}{y'_k}\right) = 1 - 0 = 1$$

The iterative formula will, near the root, approach the following expression;

$$\lim_{k \to \infty}(x_{k+1}) = \left(1 - \frac{1}{1}\right) x_{k-1} + \left(\frac{1}{1}\right) x_k = (0) x_{k-1} + (1) x_k$$

The iteration, therefore, moves away from $x_{k-1}$ with a weighing factor approaching zero and gives weighing factor of 1 to $x_k$.

Near points where the derivative of the function may approach zero, the r value approaches infinity and the root estimate weighs heavily in favor of $x_{k-1}$ rather than $x_k$. This can be seen from the following limit evaluation of r.

$$\lim_{dy/dx \to 0}(r) = 1 - \left(\frac{y_k}{y_{k-1}}\right)\left(\frac{\frac{y_k - y_{k-1}}{x_k - x_{k-1}}}{0}\right) = 1 \pm \infty = \pm\infty$$

$$\lim_{dy/dx \to 0}(x_{k+1}) = \left(1 - \frac{1}{\infty}\right)x_{k-1} + \left(\frac{1}{\infty}\right)x_k = (1)\,x_{k-1} + (0)\,x_k$$

The iteration, therefore, moves away from $x_k$ with a weighing factor approaching zero and gives weighing factor of near to 1 to $x_{k-1}$

Near points where the derivative of the function may approach zero, Newton method typically displays oscillation or off-shooting behavior. The proposed two-point modified Newton procedure however shows stability without displaying oscillation and off-shooting tendencies. For example, for the function $y = \sin(x)$, the derivative of y (i.e., $dy/dx = \cos(x)$) is zero at $x = \pi/2$. Starting Newton method near this point, say $x = (\pi/2 + 0.01)$, will result in off-shooting of the iteration to a different domain.

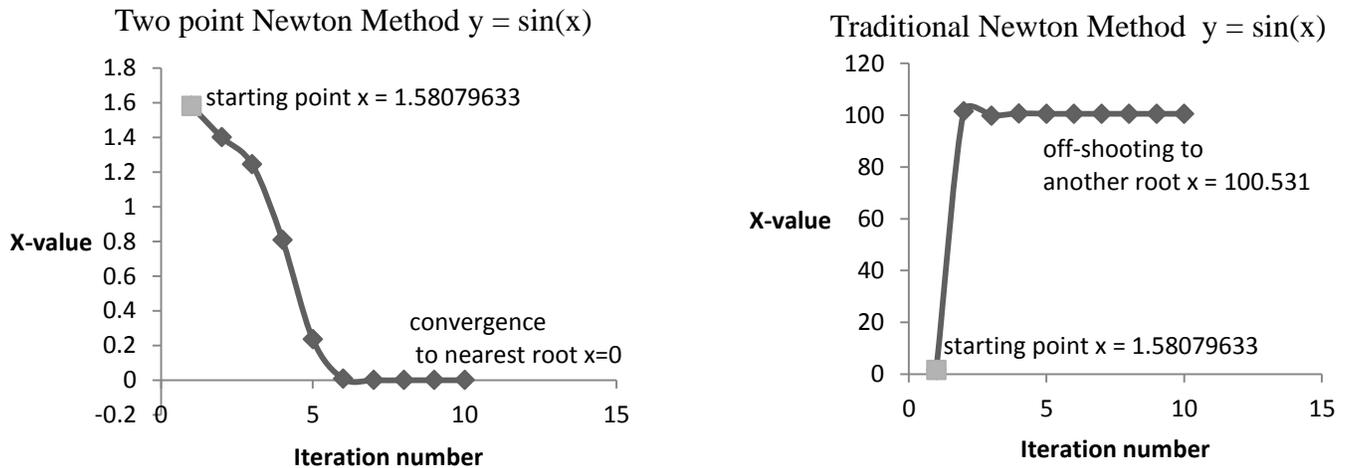

Figure 2. Stability of the proposed Two-Point Newton method near points of zero derivatives.

As shown in Figure 2, for a starting point x= 1.58079633 where the function y= sin(x) has near zero derivative, the traditional Newton method off-shoots to a root further from the nearest root of x=0. The newton method as such converges to x=100.531. On the other hand, the proposed two point Newton method shows stability and converges to the nearest root of x=0.

Another example of the stability of the proposed method is where Newton method diverges instead converging to the desired root. Figure 3 and Figure 4 show comparison of the iterative values of x for the proposed method and the traditional Newton method for the function y = x$^{1/3}$ which is often used to illustrate the pathological condition with respect to application of Newton's method which leads to an oscillating divergence to infinity. By contrast the proposed method displays an oscillating convergence to the desired root albeit with slower rate of convergence of order 1.0.

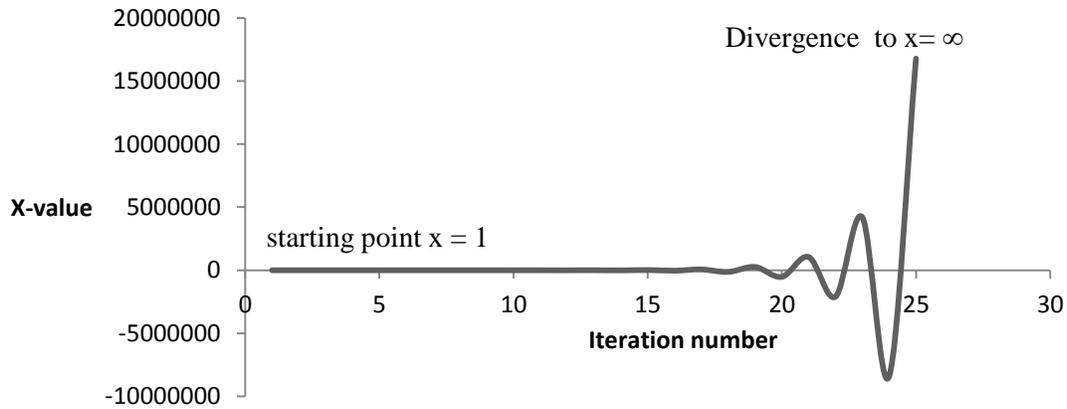

Figure 3. Application of Newton method to y= x$^{1/3}$ leading to oscillating divergence

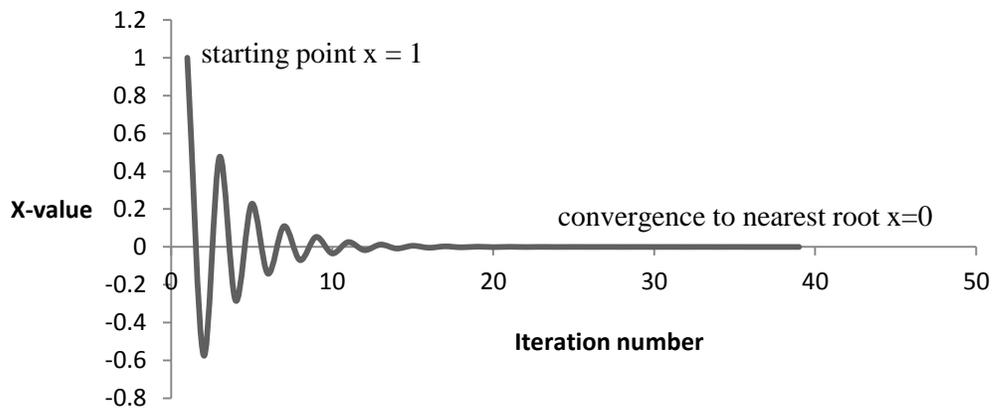

Figure 3. Application of the proposed Two–Point Newton method to y= x$^{1/3}$ leading to oscillating convergence.

## 5. Application Examples

Equations used to test efficiency of root finding methods are used here to evaluate the number of iterations required to reach to a specified level of convergence (Table 1). The stopping criterion used for the iteration process is given by:

$$|x_k - x_{k-1}| + |y_k| < 10^{-15} \tag{6}$$

The rate of convergence towards the root x = r for each step of the iteration is evaluated using the formula:

$$C_k = \frac{Log\ (|E_{k+1}|)}{Log\ (|E_k|)} = \frac{Log\ (|X_{k+1} - r|)}{Log\ (|X_k - r|)} \tag{7}$$

Table 1 below shows comparison of the proposed two-point Newton method with those of Newton and secant methods for a number of equations used to test efficiency of root finding methods elsewhere. A super-quadratic convergence with which the proposed method converges to the root is mostly evident with a $C_k$ value being close to 2.414 during most of the iterations. It can also be seen from Table 1 that less number of iterations are required to reach convergence for the proposed method than those required for Newton and secant methods. For example, for equations with multiple roots such as $y = (x - 1)^6 - 1$, shown in Table 1, Newton method displays linear convergence while the proposed method converges super-quadratically with an order of 2.414. Similar order of convergence is also observed for the equation: $y = e^{x^2+7x-30} - 1$.

### 5.1 Examples where the proposed method works while Newton method fails

The advantage of the use the proposed Two-Point Newton method is shown for cases where the Newton method and in several cases also the secant method fail to converge to the root. Table 2 below shows the results of the iteration for several examples of equations. In all of the examples listed, Newton method fails to converge whereas the proposed two-point Newton method converges. Secant method also fails to converge in several of the examples cited in Table 2. For example, in the case of y = Log(x), Newton method starting with x=3 fails with the second iteration because the estimated x value is a negative number whose logarithm is undefined. For the case y= arctan(x) Newton method always diverges to increasingly large number for any starting x value while the proposed two-point Newton method converges super qudratically near the root. For y= $x^{1/3}$, Newton method also diverges whereby each iteration gives an estimate of the root which is twice the previous value and with alternate signs leading to oscillating divergence. On the other hand application of the proposed two-point modification of Newton method leads to oscillating convergence at somehow reduced rate of convergence (of order 1 as explained above).

For the equation $y = 10xe^{-x^2} - 1$, both Newton and secant methods continue to diverge to increasingly large x values while the proposed two-point Newton method shows stability of convergence. For the polynomials of different degrees cited in Table 2, Newton and secant methods display oscillation for the starting values shown in the Table 2 while the proposed Two-Point Newton method displays stability and super-quadratic convergence for iteration near the root.

Table 1. Comparison of result of iterations of the Two-Point Newton method with Newton and Secant Methods.

| Function | Root | Starting point | Secant Method | Newton Method | Two-point Newton Method |
|---|---|---|---|---|---|
| $y = [\sin(x)]^2 - x^2 + 1$ | -1.404491648215340 |  | 10 | 8 | 6 |
|  |  | 1 | 9 | 6 | 5 |
|  |  | -1 | 9 | 7 | 5 |
|  |  | -3 | 9 | 7 | 6 |
| $y = (x-2)(x+2)^4$ | -2.0000000000000 | -3 | 168 | 119 | 83 |
|  |  | 1.4 | 116 | 81 | 60 |
| $y = (x-1)^6 - 1$ | 2.00000000000000 | 1.5 | 25 | 17 | 8 |
|  |  | 2.5 | 12 | 8 | 6 |
|  |  | 3.5 | 16 | 11 | 8 |
| $y = \sin(x) \cdot e^x + \ln(x^2+1)$ | -0.603231971557215 | -0.8 | 8 | 7 | 5 |
|  |  | -0.65 | 8 | 5 | 4 |
| $y = e^{x^2+7x-30} - 1$ | 3.000000000000000 | 4 | 29 | 20 | 14 |
|  |  | 4.5 | 39 | 28 | 18 |
| $y = x - 3\ln(x)$ | 1.857183860207840 | 2 | 8 | 5 | 4 |
|  |  | 0.5 | 11 | 8 | 5 |

.

Table 2. Results of application of the proposed method for cases in which Newton or Secant method fail to converge.

| Function | Root | Comparison of number of iterations required | | | |
|---|---|---|---|---|---|
| | | Starting point | Secant Method | Newton Method | Two-point Newton Method |
| $y = -x^4 + 3x^2 + 2$ | 1.887207676120680 | 1 | 11 | Oscillates | 7 |
| | | 0.5 | 23 | Oscillates | 6 |
| $y = \log(x)$ | 1.000000000000000 | 3.0 | Fails | Fails | 5 |
| $y = Arctan(x)$ | 0.0000000000000 | 3.0 | Diverges | Diverges | 6 |
| | | -3.0 | Diverges | Diverges | 6 |
| $y = x^5 - x + 1$ | -1.167303978261420 | 2.0 | Oscillates | Oscillates | 12 |
| | | 3.0 | 14 | Oscillates | 15 |
| $y = 0.5x^3 - 6x^2 + 21.5x - 22$ | 4.00000000000000 | 3.0 | 10 | Oscillates | 7 |
| | | 5.0 | 8 | Oscillates | 6 |
| $y = x^{1/3}$ | 0.00000000000000 | 1.0 | Oscillates | Diverges | 101 |
| | | -1.0 | Oscillates | Diverges | 101 |
| $y = 10xe^{-x^2} - 1$ | 1.679630610428450 | 3.0 | Diverges | Diverges | 8 |
| | 0.101025848315685 | -1.0 | Diverges | Diverges | 11 |

## 6. Conclusion

A numerical procedure of root finding using two-point modification of Newton method has been presented. It is proved that the method has a super-quadratic convergence of order about 2.414. The method is based on application of Newton iteration formula by taking as the independent variable the cotangent of the angle between the line connecting the two successive points of iteration with the vertical and as the dependent variable the given function y = f(x). The resulting iteration formula for root estimation is shown to be the weighted sum of the estimates of the two previous iterations with a weighing factor that penalizes the iteration point having undesirable characteristics such as a near zero derivative. For example, near a point where the derivative is zero, the weighing factor for that point will be near zero effectively moving the iteration away from that undesirable point.

Application examples have been given that demonstrated that the proposed method requires fewer number of iterations for convergence to a root than the traditional Newton and secant methods. The method offers particular advantage for cases where the traditional Newton method and its variants of various order convergence may not converge. A number of examples are given where the proposed method converges to a root in a stable manner without oscillation, divergence or off-shooting from the desired domain of the root whereas the Newton method, and in several of the cases also the secant method, display undesirable behaviors which prevent convergence to the desired root.

In terms of computational effort, the proposed method requires one function evaluation and one derivative evaluation at each step of the iteration except for the first step where two functional evaluations are required. The proposed method, therefore, requires essentially the same number of functional evaluations as the traditional Newton method while offering a super-quadratic convergence of order 2.414.